\title{\bf{Singularities in Positive Characteristic}}
\author{
       \bf{ Gert-Martin Greuel }\\[1.0ex]{\em Dedicated to Antonio Campillo on the occasion of his 65th birthday}        
           }
\date{\today}

\documentclass[11pt]{article}
\usepackage{amsmath, amsthm, amssymb, mathrsfs, amscd, amsthm, amscd,amsfonts,graphicx}
\usepackage{indentfirst,url,xypic}
\usepackage{lipsum} 
\usepackage[all]{xy}
\usepackage{url}
\usepackage[colorlinks,plainpages]{hyperref}
\usepackage{verbatim}
\usepackage{fancyvrb}
\usepackage{listings}
\hypersetup{
	colorlinks=true,
	linkcolor=blue,
	filecolor=magenta,      
	urlcolor=cyan,
}

\usepackage{epstopdf}
\usepackage{pgfplots}
\usepackage{tikz}
\usetikzlibrary{arrows,matrix,fit,positioning}

\DeclareMathOperator{\Spec}{Spec}
\DeclareMathOperator{\Ker}{ker}

\DeclareMathOperator{\Def}{Def}
\DeclareMathOperator{\Aut}{Aut}

\newtheorem{Definition}{ Definition}[section]
\newtheorem{Theorem}[Definition]{Theorem}

\newtheorem{Proposition}[Definition]{Proposition }
\newtheorem{Corollary}[Definition]{Corollary}
\newtheorem{Example}[Definition]{Example}

\newcommand{\Q}{\mathbb{Q}}

\newcommand{\C}{\mathbb{C}}
\newcommand{\Z}{\mathbb{Z}}

\newcommand{\kk}{{\mathcal K}}

\newcommand{\ko}{{\mathcal O}}

\newcommand{\kr}{{\mathcal R}}

\begin{document}
\maketitle

\begin{abstract}
In this survey paper we give an overview on some aspects of singularities of algebraic varieties over an algebraically closed field of arbitrary characteristic. We review in particular results on equisingularity of plane curve singularities, classification of hypersurface singularities and determinacy  of arbitrary singularities. The  section on equisingularity has its roots in two important early papers by Antonio Campillo. One emphasis is on the differences between positive and zero characteristic and on open problems.
\end{abstract}




\renewcommand{\contentsname}{Table of Contents}
\tableofcontents

\section{Historical Review}\label{review}

Singularitiy theory means in this paper the study of systems of polynomial or analytic or differentiable equations locally at points where the Jacobian matrix has not maximal rank. This means that the zero set of the equations at these points is not smooth. The points where this happens are called  {\em singularities} of the variety defined by the equations. Singularities have been studied since the beginning of algebraic geometry, but the establishment of their own discipline arose about 50 years ago.\\

Singularity theory started with fundamental work of {\em Heisuke Hironaka} on the resolution of singularities (1964),  
{\em  Oskar Zariski's} studies in equisingularity, (1965-1968), 
{\em Michael Artin's} paper on isolated rational singularities of surfaces (1966), and 
the work by {\em Ren\'e Thom, Bernard Malgrange, John Mather,...} on singularities of differentiable mappings.

It culminated in the 1970ties and 1980ties with the work of
{\em John Milnor}, who intorduced what is now called the Milnor fibration and the Milnor number (1968), 
{\em Egbert Brieskorn's} discovery of exotic spheres as neighborhood boundaries of isolated hypersurface singularities (1966) and the connection to Lie groups (1971),
{\em Vladimir Arnold's} classification of (simple) singularities (1973), and
many others, e.g. {\em Andrei Gabrielov, Sabir Gusein-Zade, Ignaciao Luengo, Antonio Campillo, C.T.C. Wall, Johnatan Wahl, L\^e D\~ung Tr\'ang, Bernard Teissier, Dierk Siersma, Joseph Steenbrink, ...}.

Besides the work of Artin, this was all in characteristic 0, mostly even for convergent power series over the complex or real numbers.\\

The first to study systematically "equisingular families" over a field of positive characteristic was {\em Antonio Campillo} in his thesis, published as Springer Lecture Notes in 1980. 

\section{Equisingularity}\label{equisingularity}

In the 1960's O. Zariski introduced the concept of equisingularity in order to study the resolution of hypersurface singularities by induction on the dimension. His idea can roughly described as follows:

\begin{itemize}

\item To resolve the singularities of $X$ consider a generic projection $X\to \C$. 

\medskip
\item If the fibres are an  ''equisingular'' family, then the resolution of a single fibre should resolve the nearby fibres simultaneously and then also the total space $X$.

\item If the fibres are plane curves then  {''equisingular''} means that the combinatorics of the resolution process of the fibre singularities is constant. 

Equivalently: the Puiseux pairs of each branch and the pairwise intersection numbers of the different branches are the same for each fibre or, the topological type of the fibre singularities is constant. 

\medskip

\item Zariski's idea works if the fibres are plane curves, but not in general. Nevertheless, equisingularity has become an independent research subject since then.

\end{itemize}


\subsection{Hamburger Noether expansions}\label{subsec.HN}

Let me now describe Campillo's early contribution to equisingularity.  There are two important papers by  Antonio Campillo: 
\begin{itemize}
\item {\em Algebroid Curves in Positive Characteristic} (\cite{Ca80})
\item
{\em Hamburger--Noether expansion over rings} (\cite {Ca83}).
\end{itemize}

The first was Campillo's thesis and appeared as Springer Lecture Notes in 1980  and 
is now the standard reference in the field. The second paper is however less known but perhaps even more important.\\

For the rest of the paper let
\textbf{$K$ denote an algebraically closed field of characteristic $p$ $\geq 0$}, unless otherwise stated.\\

We consider in this section reduced plane curve singularities. Over the complex numbers these are 1-dimensional complex germs germs $C \subset \C^2$ with isolated singularity at 0, given by a convergent power series $f \in \C\{x,y\}$ with $C$ the germ of the set of zeros $V(f)$ of $f$. If $K$ is arbitrary, a plane curve singularity is given by a formal power series $f \in  K[[x,y]]$ with  $C= V(f) = \Spec K[[x,y]]/\langle f\rangle$. $C$ and $f$ are also called {\bf algebroid plane curves}. \\

A reduced and irreducible algebroid plane curve $C$ can be given in two ways:
\begin{itemize}
\item by an {\em equation} $f =0$, with $f$ irreducible in the ring $K[[x,y]]$ 
\item by a  {\em parametrization} $x=x(t), y=y(t)$ with $\langle f\rangle=\Ker(K[[x,y]] \to K[[t]])$
\end{itemize}


\noindent {\bf Case $p=0$}
\begin{itemize}
\item In this case $C$ has a special parametrization (the {\bf Puiseux expansion})
\end{itemize}
$\begin{array}{ll}
\ \ \ \ x=t^n & : n=\text{ord}(f)=\text{mult} (C)\\
\ \ \ \ y= c_mt^m+c_{m+1}t^{m+1}+\cdots & : m\geq n, c_i\in K 
\end{array}$\\

\medskip
Here ord($f$) is the {\bf order} of $f$, also denoted the {\bf multiplicity} of $f$, i.e. lowest degree of a non-vanishing term of the power series $f$.

If $f = f_1 \cdot ... \cdot f_r$  is reducible (but reduced) with $f_i$ irreducible, we consider the parametrization of each {\bf branch} $f_i$ of $f$ individually.
\bigskip

The Puiseux expansion determines the \textbf{characteristic exponents} of $C$
(equivalently the \textbf{Puiseux pairs} of $C$ ). 

\begin{itemize}
\item These data is called the \textbf{equisingularity type}  ({\bf es--type}) of the irreducible $C$.

\item For a reducible curve $C$ the {\bf es--type} is defined by the  {\bf es--types of the branches} and the  {\bf pairwise intersection multiplicities} of different branches.
\item Equivalently by the {\bf system of multiplicities} of the reduced total transform in the resolution process of $C$ by successive blowing up points.
\item Two curves with the same es--type are called  \textbf{equisingular}

\end{itemize}

For the case $K=\C$ and  $f, g\in \C\{x,y\}$ we have the following nice topological interpretation of equisingularity:

\begin{itemize}

\item V(f) and V(g) are equisingular $\Leftrightarrow$ they are (embedded) \textbf{topologically equivalent},
i.e. there is a homeomorphism of triples $h:(B_\varepsilon, B\cap V(f), 0)\xrightarrow{\sim}(B_\varepsilon, B_\varepsilon\cap V(g), 0)$, with $B_\varepsilon \subset \C^2$ a small ball around 0 of radius $\varepsilon$ (cf. fig. 1).
\end{itemize}

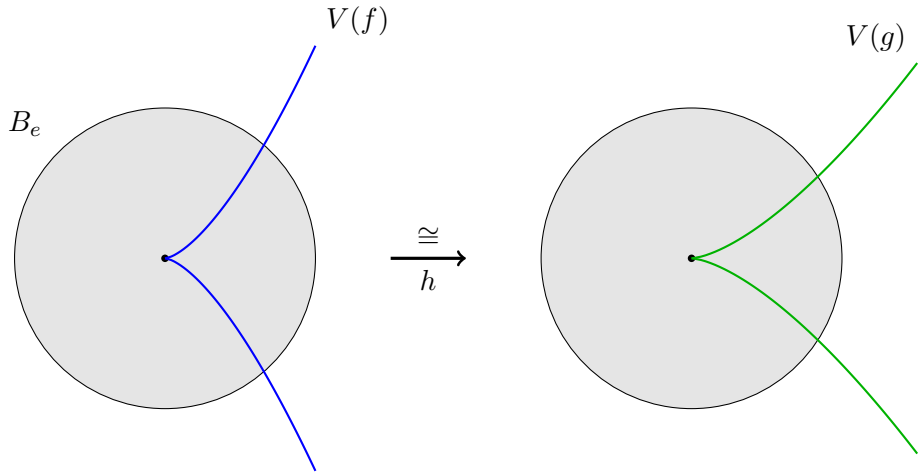
\begin{center}
\begin{tikzpicture}[scale=1.0, every node/.style={transform shape}]

\begin{scope}
\draw[fill=gray, fill opacity=0.2] (0cm,0cm) circle(2cm);
\node[above left] at (-1.5cm,1.5cm) {$B_e$};
\fill[thick] (0cm,0cm) circle(0.05cm);
\draw[thick,domain=0:2,smooth,variable=\y,draw=blue] plot
(\y,{sqrt(\y^3)}) node[above right] {$V(f)$};
\draw[thick,domain=0:2,smooth,variable=\y,draw=blue] plot
(\y,{-sqrt(\y^3)});
\end{scope}

\begin{scope}[shift={(7,0)}]
\draw[fill=gray, fill opacity=0.2] (0cm,0cm) circle(2cm);
\fill[thick] (0cm,0cm) circle(0.05cm);
\draw[thick,domain=0:3,smooth,variable=\y,draw=green!70!black] plot
(\y,{1/2*sqrt(\y^3)}) node[above left] {$V(g)$};
\draw[thick,domain=0:3,smooth,variable=\y,draw=green!70!black] plot
(\y,{-1/2*sqrt(\y^3)});
\end{scope}

\draw[->, very thick] (3,0) -- (4,0) node[above,pos=0.5]{$\cong$}
node[below,pos=0.5]{$h$};
\end{tikzpicture}
\end{center}



\noindent{\bf Case $p>0$}\\

The resolution process of $C$ by successive blowing up points exists as in the case $p=0$. There exists also a parametrization of $C$, but a \textbf{Puiseux expansion does not exist} if $p | n$, $n$ the multiplicity of $C$.

\begin{itemize}
 \item  We define the \textbf{equisingularity type (es-type)} of $C$,  by the \textbf{system of multiplicities} of the resolution as in characteristic 0. 

 \item  Instead of the Puiseux expansion another special parametrization exists and can be computed from any parametrization (or any equation) of $C$, the \textbf{Hamburg--Noether (HN) expansion}. It is  determined by a chain of relations obtained from successive divisions by series of lower order (assume $ord(x) \leq ord(y)$) as follows:
 
$
\begin{array}{lclcr}
y & = & a_{01}x \ +\cdots + a_{0h}x^h+x^h z_1, & ord(z_1) < ord(x)\\
x & = & a_{12}z_1^2 + \cdots +a_{1h_1}z_1{h_1}+z_1^{h_1}z_2 \\
z_1 & = & a_{22}z_2^2 + \cdots +a_{2h_2}z_2{h_2}+z_2^{h_2}z_3 && ({\bf HN(C)}) \\
& & \vdots & &\\
z_{r-1}& = & a_{r2}z_r^2+\cdots   \ \ \in K[[z_r]]
\end{array}
$
\end{itemize}

We do not have Puiseux pairs, but we have characteristic exponents. Campillo  defines the {\bf characteristic exponents for $C$} from HN(C). \\

By substituting backwards, we get from HN(C) a parametrization of $C$:
\begin {center}
 $x=x(z_r), \ \ y = y(z_r) \ \in K[[z_r]]$\\
\end {center}
Note that the uniformizing parameter $z_r$ is a rational function of the coordinates $x,y$. It does not involve roots of unity as the uniformizing parameter for the Puiseux expansion. Moreover, computationally the Hamburg--Noether expansion is preferred to the Puiseux expansion as it needs the least number of field extensions if one wants to compute the es-type for an algebroid curve defined over a non algebraically closed field (such as $\Q$). This is implemented in SINGULAR \cite{DGPS16}.\\

For an arbitrary algebraically closed field $K$ Campillo defines the complex model of $C$ as follows:
\begin{itemize}
\item  Let $F:K\to \C$ be any map with $F(a)=0\Leftrightarrow a=0$
 (e.g. $F(0)=0, F(a) = 1$ for $a\neq 0$).
 
A \textbf{complex model}  $C_\C$  of the curve $C$ is obtained from HN(C) by the HN-expansion
\begin {center}

HN($C_\C$): replace  $a_{ij}$  in HN($C$) by  $F(a_{ij}).$
\end {center}
\end{itemize}
\bigskip

\begin {Theorem} {\bf(Campillo, \cite {Ca80})} \label{th2.1}
\begin{enumerate}
\item The characteristic exponents of $C_\C$ do not depend on the complex model.
\item They are a complete set of invariants of the es-type of $C$.
\item They coincide with the characteristic exponents of $C_\C$ obtained from the Puiseux espansion.
\end{enumerate}
\end {Theorem}

Note that the complex model $C_\C$ of $C$ is defined over the integers if $F$ has integer values (this is important for coding theory and cryptography).\\


We come now to the second paper of Campillo ''Hamburger--Noether expansion over rings'' mentioned above. For any ring $A$ Campillo defines a
\begin{itemize}
\item  \textbf{HN--expansion HN$_A$ over $A$}. HN$_A$ is similar to HN, but with $a_{ij} \in A$ and certain properties. 
It may be considered as a family over Spec($A$) of parametrized curves with constant es-type. If $A$ is the field $K$ then HN$_K$ coincides with HN defined above. 
\end{itemize}

If $A$ is a local $K$--algebra with maximal ideal $\mathfrak {m}_A$ 
and $A/\mathfrak{m}_A = K $, then we may take 
residue classes of the HN-coefficients $a_{ij}$ modulo $\mathfrak{m}_A$ and thus the HN-expansion over $A$ induces a deformation 
\[
X=X(z_r),\ \ Y=Y(z_r) \ \in A[[z_r]]
\]
of the parametrization over Spec$(A)$ 
\[
x=x(z_r), \ \ y=y(z_r)\in K[[z_r]], \ \text{with } x, y =  X ,Y  \text {mod }  \frak{m} _A
\]
of an irreducible plane algebroid curve $C$. 

If A is irreducible, the es--types of the curve $C$ parametrized by $x(z_r), y(z_r)$ over $K$ and of the curve defined by the parametrization $X(z_r), Y(z_r)$ over the quotient field Quot$(A)$ coincide. \\

We have the following important theorem, saying that for a fixed equisingularity type there exists some, in a sense ''totally versal", equisingular family $X \to Y$ of $\Z$--schemes such that for any field $K$ the following holds: any equisingular family of algebroid curves over $K$ can be induced from $X \to Y$. More precisely, Campillo proves:\\

\begin {Theorem} {\bf(Campillo, \cite {Ca83})}\label{th2.2}

\begin{enumerate}
 \item [(1)] Let $C$ be irreducible and $E$ the es--type of $C$. Then there exists a morphism of $\Z$--schemes $\pi:X\to Y$ with section $\sigma:Y\to X$ s.t. for any algebraically closed field $K$ the base change $\Z\to K$ induces a family $X^K\to Y^K$ with section $\sigma^K$ such that the following holds:
\begin{enumerate}
\item [(i)] $Y^K$ is a {\em smooth, irreducible} affine algebraic variety over $K.$
\item [(ii)] For any closed point $y\in Y^K$ the induced family 
\[
\text{Spec}(\widehat{\ko}_{X^K, \sigma^K(y)})\to \text{Spec}(\widehat{\ko}_{Y^K, y})
\]
is a {\bf total es--versal HN--expansion}, i.e.:

for any algebroid curve $C'$ with es$(C') =E$ and any local $K$--algebra $A$ s.t.
$HN_A$ induces $HN(C')$ mod $\frak{m}_A$, there exists a morphism $\varphi: \widehat{\ko}_{Y^K, y} \to A$ such that $HN_A$ is induced from $HN_{\widehat{\ko}_{Y^K, y}}$ via $\varphi$.

\end{enumerate}
\item [(2)] For a reducible curve $C$ the construction is extended to finite sets of HN--expansions over $A$ and the statement of (1) continues to hold.
\end{enumerate}
\end {Theorem}


\subsection{Equisingularity strata}\label{subsec.ES}

A Hamburger-Noether expansion over a ring $A$ induces an {\bf equisingular deformation of the parametrization} of an irreducible curve singularity $C$. Such a deformation of the parametrization induces a deformation of the equation as follows: 

\begin{itemize}
\item
Let  $X=X(z_r), Y=Y(z_r)$ be a HN--expansion over $A$. It is a deformation of the parametrization  $x=x(z_r), y= y(z_r)$  of an algebroid curve  $C=\text{Spec}(K[[x, y]]/\langle f \rangle)$.   

By elimination of $z_r$ from $x-X(z_r)$ and $y-Y(z_r)$, we get a power series $F\in A[[x,y]]$ with $f=F \mod \frak{m}_A$. $F$ is a deformation of the curve $C=V(f)$ (in the usual sense) over Spec$(A)$, also called a deformation of the equation.  Since it is induced from an equisingular deformation of the parametrization, we call it an equisingular deformation of the equation. Deformations are a category and the construction is functorial (cf. \cite {GLS07}).

\item In this way we get for any algebroid curve $C$ a functor
\[
\begin{array}{ll}
\chi_{es}:  & \text{equisingular--deformations of the parametrization of C} \\
& \rightarrow  \text{(usual) deformations of the equation of C}.\\
\end{array}
\]
We call the image  of $\chi_{es}$ (full subcategory) {\bf equisingular--deformations of the equation} or just  {\bf es--deformations} of $C$ over Spec($A$).
\item The construction can be generalized to reducible $C$ and a set of HN-expansions over $A$ of the branches of $C$ (with certain properties). 
\end{itemize}

More generally, any deformation $\Phi: X \to T$, equisingular or not, of the parametrization of a plane curve singularity $C$ induces a deformation of the equation by eliminating the uniformising variable. \\

{\bf Question:} Does the base space $T$ of an arbitrary deformation of $C$ admit a unique maximal subspace over which the deformation is equisingular? In other words, does there exist a unique {\bf equisingularity stratum} of $\Phi$ in $T$? \\

The answer is well--known for $K = \C$. Recall that for any $K$ and $f \in K[[x,y]]$
\begin{itemize}
\item
$ \mu(f):=\dim_K K[[x,y]]/\langle f_x, f_y\rangle$ is the {\bf Milnor number} of $f$. 
\end{itemize}

If $p$ = char($K$) = 0 $ \mu(f)$ depends only on the ideal $\langle f\rangle$ and not on the choosen generator and then we write also $ \mu(C)$ for $C=\text{Spec } K[[x, y]]/\langle f \rangle$.\\

For $K = \C$ or, more generally, if char($K$) = 0, the equisingularity stratum of any deformation $\Phi: X \to T$ of $C$ exists and is the {\bf $\mu$--constant stratum} of  $\Phi$, i.e. the set of points $t \in T$ such that Milnor number of the fibres $X_t$ is constant along some section $\sigma :T \to X$ of $\Phi$.
If $\Phi: X \to T$ is the semiuniversal deformation of $C$, the  $\mu$--constant stratum is denoted by $\Delta_{\mu}$,  
\begin{itemize}
\item $\Delta_{\mu}$ = \{$t \in T$ : $\exists$ a section $\sigma :T \to X$ of $\Phi$ s.t. $\mu(X_t, \sigma(t)) = \mu(C)$\}. 
\end{itemize}

 The restriction of $\Phi$  to $\Delta_{\mu}$ may be considered as the {\bf semiuniversal es--deformation} of $C$ in the sense that any es--deformation of $C$ over some base space $T$ can be induced by a morphism $\varphi: T \to \Delta_{\mu}$, such that the tangent map of $\varphi $ is unique. Moreover $\Delta_{\mu}$ is known to be smooth (c.f. e.g. \cite{GLS07}).\\

If $K$ has positive characteristic, the situation is more complicated. A semiuniversal es--deformation of $C$ exists, but an equisingularity stratum does not always exist. The situation is described in the following theorem and in the next subsection.

\begin {Theorem}  {\bf (Campillo, Greuel, Lossen;  \cite{CGL07})} \label{th2.3}

\begin{enumerate}
\item [(1)] The functor $\chi_{es}$ defined above is smooth (unobstructed).
\item [(2)] The functor $\underline{\Def}^{es}_C$  of (isomorphism classes) of es--deformations of the equation has a semiuniversal deformation with smooth base space $B_C^{es}$.
\item [(3)] If $\text{char}(K)=0$ then  $B_C^{es}$ 
coincides with the $\mu$--constant stratum $\Delta_{\mu}$ in the base space of the (usual) semiuniversal deformation of $C$.

\item [(4)] 
In \textbf{good characteristic} (i.e. either $p=0$ or $p>0$ does not divide the multiplicity of any branch of $C$)
there exists for any deformation of $C$ over some $T$ a unique maximal  equisingularity stratum $T_{es} \subset T$  (generalizing the $\mu$--constant stratum). 
\end{enumerate}
\end {Theorem}


\subsection{Pathologies and open problems}\label{subsec.POP}

If the characteristic is bad, i.e. $p >0$ divides the multiplicity of some branch of $C$, we have the following {\bf pathologies}:

\begin{enumerate}
\item [(1)] There eixst deformations of $C$ wich are not equisingular but become equisingular after a finite base change. We call these deformations \textbf{weakly equisingular}.

\item [(2)] Let $\Phi_C: X_C\to B_C$ denote the semiuniversal deformation of $C$. In general no unique es--stratum in $B_C$ exists. For example, let $p=2$ and $f=x^4+y^6+y^7$. Then the following holds.

There exist infinitely many smooth subgerms $B_\alpha\subset B_C$ s.t.
\begin{itemize}
\item $B_\alpha\cong B_C^{es}$.
\item All $B_\alpha$ have the same tangent space.
\item The restricton of $\Phi_C$ to $B_\alpha$ is equisingular for any $\alpha$.
\item The restriction of $\Phi_C$ to $B_{\alpha_1}\cup B_{\alpha_2}$ is not equisingular if $\alpha_1\neq \alpha_2$. Hence, a unique maximal equisingularity stratum of $\Phi_C$ does not exist.
\end{itemize}

\item [(3)] Let $B_C^{wes}$  denote the Zariski--closure of the union of all subgerms $B'\subset B_C$ with the property that $\Phi_C|{B'}$ is equisingular. Then $\Phi_C|B_C^{wes}$ is  the  \textbf {semiuniversal weakly equisingular deformation} of $C$, i.e it has the versality property for weakly es--deformations (as above for es--deformations). We have
\begin{itemize}
\item $B_C^{wes}$ is irreducible but in general not smooth
\item $B_C^{wes}$ becomes smooth after a finite (purely inseparable) base change.
\end{itemize}
\end{enumerate}

For a proof of these facts see \cite{CGL07}.
\bigskip

Let us mention the following {\bf Open Problem}:

\begin{itemize}
\item We know from Theorem \ref{th2.3}  that "$p$ = good" is a sufficient condition for $B_C^{es}= B_C^{wes}$ (and hence that $B_C^{es}$  is smooth) but it is not a necessary condition.
The problem is to find necessary conditions for $B_C^{es}=B_C^{wes}$ (these do not only depend on $p$).
\end{itemize}

We do not yet fully understand the relation between weak and strong equisingularity (even between weak and strong triviality).

\section{Classification of Singularities}\label{classification}

In this section we consider hypersurface singularities $f \in K[[x]]:=K[[x_1, \ldots, x_n]]$, again with $K$ algebraically closed and char$(K) = p\geq 0$. We recall the classical results for the classification of singularities in characteristic zero and present some more recent results in positive characteristic.

The two most important equivalence relations for power series are right equivalence and contact equivalence. The two equivalence relations lead of course to different classification results. It turns out that the classification of so called ''simple singularities''  w.r.t. contact equivalence is rather similar for $p=0$ and for $p > 0$. However, for right equivalence the classification of simple singularities in positive characteristic is surprisingly different from that in characteristic zero.\\

Let $f,g  \in K[[x]]$. Recall that 
\begin{itemize}
\item $f$ is \textbf{right equivalent} to $g  \ (f\overset{r }{\sim} g):
 \Leftrightarrow \exists\  \Phi \in \Aut(K[[x]])$ such that  $f= \Phi(g)$, 
i.e. $f$ and $g$ differ by an analytic coordinate change $\Phi$ with $\Phi(x_i)=\Phi_i $ and $f=g(\Phi_1, \cdots, \Phi_n).$

\item $f$ is \textbf{contact equivalent} to $g \ (f \overset{c}{\sim}\ g):
\Leftrightarrow \exists\ \Phi \in \Aut (K[[x]])$ and $u \in K[[x]]^\ast$ such that $f=u\cdot\Phi(g)$, 
i.e. the analytic $K$--algebras $K[[x]]/\langle f\rangle $ and $K[[x]]/\langle g\rangle$ are isomorphic.

\item $f$ is called \textbf{right--simple} (resp. \textbf{contact--simple}): 
$\Leftrightarrow \exists$ finite set $\{g_1, \ldots, g_l\}\subset K[[x]]$ of power series such that for any deformation of $f$,
\begin {center}
$F_t(x)=F(x,t)=f(x)+\sum\limits^k_{i=1} t_ih_i(x,t), $ 
\end {center}

there exists a neighbouhood $U=U(0)\subset K^k$ such that $\forall\ t\in U\ \exists\ 1\leq j \leq l$ with $F_t\overset{r }{\sim} g_j$ (resp. $F_t\ \overset{c }{\sim} \ g_j$); in other words,
$f$ has  \textbf{no moduli} or $f$ is of \textbf{finite deformation type}.
\end{itemize}


\subsection {Classification in characteristic zero}

The most important classification result for hypersurface singularities in characteristic zero is the following result by V. Arnold:

\begin{Theorem} {\bf (Arnold; \cite{AGV85})}
Let $f\in \C\{x_1, \ldots, x_n\}.$
Then $f$ is right--simple $\Leftrightarrow f$ is right equivalent to an ADE singularity from the following list:
\\

$
\begin{array}{ll}
A_k: & x_1^{k+1} + x_2^2 \ \ \ \ \ \   +q,\ \ \ \ \ k\geq 1,\quad q:=x_3^2+\cdots + x_n^2\\
D_k: & x_2(x_1^2+x_2^{k-2}) \ +q, \ \ \ \ \ k\geq 4\\
E_6: & x_1^3+x_2^4 \ \ \ \ \ \ \ \ \   +q\\
E_7: & x_1(x_1^2+x_2^3) \ \ \ \ +q\\
E_8: & x_1^3+x_2^5 \ \ \ \ \ \ \ \ \  +q
\end{array}
$
\end {Theorem}

\noindent Later it was proved that $f \in \C\{x_1, \ldots, x_n\}$ is right--simple $\Leftrightarrow f$ is contact--simple. \\

Arnold's classification has numerous applications. The list of simple or ADE singularities appears in many other contexts of mathematics and is obtained also by classifications using a completely different equivalence relation (cf. \cite {Du79}, \cite{Gr92}). One such classification result is the following.

\begin{Theorem} {\bf (Buchweitz, Greuel, Schreyer; Kn\"orrer, \cite{BGS87}, \cite{Kn87})}

$f \in \C\{x_1, \ldots, x_r\}$ is simple $\Leftrightarrow$ $f$ is of \textbf{finite CM--type}, (i.e. there are only finitely many isomorphism classes of maximal indecomposable Cohen--Macaulay modules over the local ring $\C\{x\}/\langle f\rangle$).
\end {Theorem}

\subsection {Classification in positive characteristic}
The classification of hypersurface singularities in positive characteristic started with the following result. Let $K$ be an algebraically closed field of characteristic $p>0$. 

\begin{Theorem} {\bf (Greuel, Kr\"oning \cite{GK90})}

The following are equivalent for $f \in K[[x_1, \cdots, x_n]]$:
\begin{enumerate}
\item [(1)] $f$ is contact--simple,
\item [(2)] $f$ is of finite $CM$--type.
\item [(3)] $f$ is an ADE--singularity, i.e. $f$ is contact equivalent to a power series form Arnold's list,  but with few extra normal forms in small characteristic ($p\leq5$).\\
E.g. for $p=3 \ \exists$ two normal forms for $E_6$: $E_6^0=x^3+y^4$ and $E_6^1=x^3+x^2y^2+y^5$.
\end{enumerate}
\end{Theorem}

The classification of right--simple singularities in positive characteristic remained open for many years. It turned out that the result differs substantially from the characteristic zero case. For example, in  characteristic zero there exist two infinite series $A_k$ and $D_k$ of right--simple singularities, but for any $p > 0$ there are only finitely many right--simple singularities.

\begin{Theorem} {\bf (Greuel, Nguyen \cite{GN16})}

$f \in K[[x_1, \cdots, x_n]]$ is right--simple $\Leftrightarrow$ $f$ is right equivalent to one of the following:
\begin{enumerate}
\item [(1)] $\mathbf{n=1}: \ \ \ \ \ \ A_k: x^{k+1} \ \ \qquad\qquad 1\leq k < p-1$
\item [(2)] $\mathbf{n=2, p>2}$:\\
\hspace*{2cm} $A_k: x^{k+1}+y^2\ ,\qquad 1\leq k < p-1$ \\
\hspace*{2cm} $D_k: x(y^2+x^{k-2}), \quad 4\leq k<p$\\
\hspace*{2.1cm}$E_6: x^3+y^4, \ \ \quad\qquad p>3$\\
\hspace*{2.1cm}$E_7: x(x^2+y^3), \ \qquad p>3$\\
\hspace*{2.1cm}$E_8: x^3+y^5 \ \ \ \ \quad\qquad p>5$\\
\item [(3)] $\mathbf{n>2, p>2}$\\
\hspace*{2cm} $g(x_1, x_2)+x_3^2+\cdots +x_n^2$ with $g$ from (2)\\
\item [(4)] $\mathbf{n\geq 2, p=2}$\\
\hspace*{2cm} $A_1: x_1x_2+x_3x_4+\cdots +x_{n-1}x_n\ ,\ n$ even
\end{enumerate}
\end{Theorem}

Note that Arnold proved a ''singularity determinator'' and accomplished the complete classification of unimodal and bimodul hypersurfaces w.r.t. right equivalence in characteristic zero with tables of normal forms (c.f. \cite {AGV85}). Such a singularity determinator and a classification of unimodal and bimodul  singularities in positive characteristic was achieved by Nguyen Hong Duc in \cite {Ng17}.

\subsection{Pathologies and a conjecture}

We comment on some differences of the classification in positive and zero characteristic and propose a conjecture.\\

For $f \in K[[x_1, \cdots, x_n]]$ the {\bf Milnor number} of $f$ is
\[
 \mu(f)=\dim_K K[[x]]/\langle f_{x_1}, \ldots, f_{x_n}\rangle.
 \]

If $f\in\frak{m}^2, \frak{m} = \langle x_1 \ldots, x_n\rangle$, and  $\mu(f) <\infty$, we have seen the following pathologies:

\begin{itemize}
\item For all $p>0$ there exist only finitely many right--simple singularities, in particular $\mu(f)$ is bounded by a function of $p$.
\item For $p=2$ and $n$ odd there exists no right--simple singularity.
\end{itemize}

The reason why there are so few right--simple singularities in characteristic $p>0$ can be seen from the following example: the dimension of the group $Aut(K[[x]])$ and hence the right orbit of $f$ is too small (due to $(x+y)^p = x^p + y^p$). \\

{\bf Example:} We show that $f=x^p+x^{p+1}$ is not right--simple in characteristic $p$
by showing that if $f_t=x^p+tx^{p+1} $
is right equivalent to $f_{t'}=x^p+t'x^{p+1}$ then $t=t'$. This shows that $f$ can be deformed into infinitely many different normal forms. 

To see this, assume $\Phi(f_t)=f_{t'}$ for $\Phi(x)=u_1x+u_2x^2+\cdots \in\Aut(K[[x]])$. Then
$(u_1x+u_2x^2+\cdots)^p+t(u_1x+u_2x^2+\cdots)^p(u_1x+\cdots)=x^p+t'x^{p+1}$ and hence
$u_1^p - 1 = (u_1 - 1)^p = 0.$ This implies $u_1=1, \ tu_1^{p+1}=t'$ and $t=t'.$\\

The classification suggests the following conjecture that is in strict contradiction to the characteristic 0 case:\\

\textbf{Conjecture:}  \\
Let char$(K)=p>0$ and $f_k\in K[[x_1,\ldots,x_n]]$ a sequence of isolated singularities with Milnor number $\mu(f_k) < \infty$ going to infinity if $k\to \infty$. Then the right modality of $f_k$ goes to infinity too. 

The conjecture was proved by Nguyen H. D. for unimodal and bimodal singularities and follows from the classification. He shows in \cite{Ng17} that $\mu(f) \leq 4p$ if the right modality of $f$ is less or equal to 2.

\section{Finite determinacy and tangent image}\label{group actions}

A power series is finitely determined (for a given equivalence relation) if it is equivalent to its truncation up to some finite order. For the classification of singularities the property of being finitely determined is indispensable. In this section we give a survey on finite determinacy in characteristic $p \geq 0$, not only for power series but also for ideals and matrices of power series. We consider algebraic group actions and their tangent maps where new phenomena appear for $p>0$, which lead to interesting open problems. 

\subsection{Finite determinacy for hypersurfaces}
A power series $f \in K[[x_1, \cdots, x_n]]$  is called right (resp. contact) $\mathbf{k}$\textbf{--determined} if for all 
$g\in K[[x]]$ such that $j_k(f)=j_k(g)$ we have $f\overset{r}{\sim} g$ (resp. $f\overset{c}{\sim} g$).
$f$ is called {\bf finitely determined}, if it is k--determined for some $k < \infty$.

\begin{itemize} 
\item Here 
\[
j_k: K[[x]]\to J^{(k)}:=K[[x]]/\frak{m}^{k+1}
\]
denotes the canonical projection, called the $k$--jet. Usually we identify $j_k(f)$ with the power series expansion of $f$ up to and including order $k$.

\item   We use ord$(f)$, the {\bf order} of $f$ (the maximum $k$ with $f \in \frak{m}^k$) and the {\bf Tjurina number} of $f$,
\[
 \tau(f):=\dim_K K[[x]]/\langle f, f_{x_1}, \ldots, f_{x_n}\rangle.
 \]
\end {itemize} 


\begin{Theorem} {\bf (Boubakri, Greuel, Markwig, Pham \cite {BGM12}, \cite{GP17})}

For $f \in \frak{m}$ the following holds.
\begin{enumerate}
\item $f$ is finitely right determined $\Leftrightarrow \mu(f)< \infty$. \\
If this holds, then $f$ is right ($2\mu(f)-$ord$(f) +2$)--determined.
 \item $f$ is finitely contact determined $\Leftrightarrow \tau(f)< \infty$. \\
 If this holds, then $f$ is contact ($2\tau(f)-ordf(f) +2$)--determined.
\end{enumerate}
\end{Theorem}

If the characteristic is 0 we have better bounds: $f$ is right ($\mu(f)+1$)--determined (resp. ($\tau(f)+1$)--determined) if $\mu(f) < \infty$ (which is equivalent to $\tau(f) < \infty$ for $p=0$), see \cite{GLS07} for $K=\C$ and use the Lefschetz principle for arbitrary $K$ with $p=0$.
The proof in characteristic $p>0$ is more difficult than for $p=0$, due to a pathology of algebriac group actions in positive characteristic, which we address in section 4.3.

\begin {itemize}
\item We can express right (resp. contact) equivalence by actions of algebraic groups. We have
\[
f\overset{r }{\sim} g \ (\text{resp.} f\overset{c }{\sim} g)\Leftrightarrow f\in  \text{orbit } G\cdot g \text{ of } g \text{ with groups}
\]
\[
G=\kr:=\Aut(K[[x]]) \textbf{ right group }     (\text{for } \overset{r }{\sim}),
\]
\[
\ \ \ G=\kk:=K[[x]]^\ast\rtimes \kr \textbf{ contact group }  (\text{for } \overset{c }{\sim}),
\]
where $G$ acts as $(\Psi=\Phi, f)\mapsto \Phi(f)$ (resp. $(\Psi=(u, \Phi), f)\mapsto u\cdot \Phi(f)$).

\item G is not an algebraic group, but the $k$--jet $G^{(k)}$ is algebraic and the induced action on $k$--jets 
\[
G^{(k)}\times J^{(k)}\to J^{(k)}\ ,\ (\Psi, f)\mapsto j_k(\Psi\cdot f)
\]
is an algebraic action.
If $f$ is finitely determined then, for sufficiently large $k$,  $f$ is (right resp. contact) equivalent to $g$ iff $j_k(f)$ is in the $G^{(k)}$--orbit of $j_k(g)$.
\end {itemize}

The determination of the tangent space $T_f(G^{(k)}f)$ of the orbit of $G^{(k)}$ at $f$ is important, but there is a big difference for $p=0$ and $p >0$. 

Consider  the tangent map $T{o_k}$ to the orbit map  
\[
o_k : G^{(k)} \to G^{(k)} f  \subset K[[x]],\ \Psi\mapsto j_k(\Psi\cdot f)
\]
and denote the image of $T{o_k} : T_e(G^{(k)}) \to T_f(G^{(k)} f) $ by $\widetilde{T}_f(G^{(k)}f)$. We have 
\[
\widetilde{T}_f(G^{(k)}f) \subset T_f(G^{(k)}f)
\]
with equality if $p=0$, but strict inclusion may happen if $p<0$ as we shall see below.
$ \widetilde{T}_f(G^{(k)}f)$ and  $T_f(G^{(k)}f)$ are inverse systems and we define the inverse limits as 
\[
T_f(Gf):=\mathop {\lim }\limits_{\mathop {\longleftarrow} \limits_{k\ge 0} } T_{f}(G^{(k)}f))\subset K[[x]],
\textbf{ tangent space}, \text{ and} 
\]
\[
\widetilde T_f(Gf):=\mathop {\lim }\limits_{\mathop {\longleftarrow} \limits_{k\ge 0} } \widetilde T_{jf}(G^{(k)}f))\subset K[[x]], \textbf{ tangent image}
\]
to the orbit $Gf$ of $G$. \\

The tangent images for $G = \kr$ and $G =\kk$ can be easily identified:
\[
\widetilde{T}_f(\kr f)=\frak{m}\langle \frac{\partial f}{\partial x_1}, \ldots, \frac{\partial f}{\partial x_n}\rangle,
\]
\[
\widetilde{T}_f(\kk f)=\langle f\rangle +\frak{m} \langle \frac{\partial f}{\partial x_1}, \ldots, \frac{\partial f}{\partial x_n}\rangle.
\]


If char$(K)=0$ then the orbit map $o_k$ is separable, which implies
$T_f(Gf)=\widetilde{T}_f(Gf)$. Moreover, in any characteristic we have:
\begin {itemize}
\item If the tangent space and the tangent image to $Gf$ coincide (e.g. if char$(K)=0$), then $f$ is finitely determined if and only if
\[  \text{dim}_K K[[x]]/\widetilde{T}_f(Gf) < \infty.
\]
\end {itemize}


\subsection{Finite determinacy for ideals and matrices}
We generalize the results of the previous section to ideals and matrices. Consider matrices 
\[
A = [a_{ij}] \in  M_{r,s}: = \text{Mat}(r, s, K[[x_1, \cdots, x_n]]) \text{ with }  r \geq s, 
\]
and the group
\[
G = \text{GL}(r, K[[x]])\times \text{GL}(s, K[[x]])\rtimes \text{Aut}(K[[x]])
\]
acting on $M_{r,s}$ in the obvious way:
\[
(U, V, \Phi, A)\mapsto U\cdot \Phi(A)\cdot V=U\cdot [a_{ij}(\Phi(x))] \cdot V.
\]
If $r=s=1$ and $A=[f]$ then $GA = \kk f$ and the considerations of this section generalize contact equivalence for power series.\\

$A$  is called $G$ {\bf k--determined}  if for all $B \in M_{r,s}$ with $j_k(A)=j_k(B)$ we have $B \subset G \cdot A$.
$A$ is {\bf finitely} $G$--{\bf determined}, if $A$ is $G$ k--determined for some $k < \infty$.\\

As in the case of one power series we have:
\begin{itemize}
\item the induced action of $G$ on jet--spaces gives algebraic group actions,
\item the tangent image to the orbit of $G$ is contained in the tangent space
\[
\widetilde{T}_A(GA)\subset T_A(GA) \text{ with }  "=" \text{ if }  p= \text{char}(K)=0,
\]
\item In any characteristic $\widetilde{T}_A(GA)$ can be computed in terms of the entries of $A$ and their partials, but $T_A (GA)$ in general not if $p>0$.
\end{itemize}


\begin{Theorem} {\bf (Greuel, Pham \cite{GP16})}
\begin{enumerate}
\item [(1)] If $\dim_K(M_{r,s}/\widetilde{T}_A(GA))<\infty\Rightarrow A$ is finitely $G$--determined
 (in particular, the orbit $GA$  contains a matrix with polynomial entries).
\item [(2)] If $A$ is finitely $G$--determined $\Rightarrow \dim_K(M_{r,s}/T_A(GA))<\infty$
\end{enumerate}
\end{Theorem}

In general we do not know whether $\dim_K(M_{r,s}/\widetilde{T}_A(GA))<\infty$ is nessecary for finite $G$--determinacy of $A$ for $p>0$, except for the case of 1--column matrices:

\begin{Theorem} {\bf (Greuel, Pham \cite{GP17})}
Let $p\geq 0$. If $A\in M_{r, 1}$, then $A$ is finitely $G$--determined $\Leftrightarrow \dim_K M_{r,1}/\widetilde{T}(GA)<\infty$.
\end{Theorem}

Two ideals $I, J \subset K[[x_1, \ldots, x_s]]$ are called {\bf contact equivalent} $\Leftrightarrow K[[x]]/I\cong K[[x]]/J$.
Since $G$--equivalence for 1--column matrices is the same as contact equivalence for the ideals generated by the entries of the matrices, we have 

\begin{Corollary}
Let $I=\langle f_1, \ldots, f_r\rangle \subset \frak{m}$ an ideal with $r$ the minimal number of generators of $I$ and let $I_r$ be the ideal generated by $r\times n$--minors of the Jacobian matrix $[\frac{\partial f_i}{\partial x_j}]$.

\begin{enumerate}
\item [(1)] $r\geq n: I$  is finitely contact determined \\
\hspace*{1.5cm} $\Leftrightarrow \dim_K(K[[x]]/I)<\infty$
\item [(2)] $r\leq n: I$ is finitely contact determined\\
\hspace*{1.5cm} $\Leftrightarrow \dim_K(K[[x]]/I+I_r)<\infty$
\end{enumerate}
\end{Corollary}

\begin{Theorem} {\bf (Greuel, Pham \cite{GP17})}
For an ideal $I=\langle f_1, \ldots, f_r\rangle \subset \frak{m} K[[x_1, \ldots, x_n]]$ with $r$ the minimal number of generators of I, $ r\leq n$, the following are equivalent in any characteristic:
\begin{enumerate}
\item [(1)] $I$ is finitely contact determined,
\item [(2)] The Tjurina number $\tau(I):=\dim_K K[[x]]^r/(I\cdot K[[x]]^r+\langle\frac{\partial f}{\partial x_1}, \cdots, \frac{\partial f}{\partial x_n}\rangle)<\infty$, with $\frac{\partial f}{\partial x_i} =( \frac{\partial f_1}{\partial x_i}, \cdots,\frac{\partial f_r}{\partial x_i}) \in K[[x]]^r$,
\item [(3)] $I$ is an isolated complete intersection singularity.
\end{enumerate}
\end{Theorem}

For the proof of (1) $\Rightarrow$ (2) we need a result about Fitting ideals, which is of independent interest.

\begin{Proposition} {\bf (Greuel, Pham \cite{GP17})}
Let $A \in M_{r,s}$ be finitely $G$--determined and let $I_t\subset K[[x_1, \cdots, x_n]]$ be the Fitting ideal generated by the $t\times t$ minors of $A$. Then $I_t$ has maximal height, i.e.
\[
ht(I_t)=\min\{s, (r-t+1)(s-t+1)\}, \ t = 1, \cdots, n.
\]
\end{Proposition} 

\subsection{Pathology and a problem}

We show that  $\widetilde{T}_f(Gf)\subsetneqq T_f(Gf)$ may happen in positive characteristic:

\begin{Example} 
Let $ G=\kk\ ,\ f=x^3+y^4\ ,\ \text{char}(K)=3$.
We compute (using SINGULAR, see \cite {GP17a}):
\begin{itemize}
\item $f$ is contact 5--determined
\item $\dim_K\widetilde{T}_f(\kk^{(5)}f)=11$
\item $\dim_K T_f(\kk^{(5)}f)=12$
\end{itemize}
\end{Example}

For the computation of $\widetilde {T}_f $ we use the formula from section 4.1 but since the tangent space $T_f$ has no description in terms of $f$ and $\frac{\partial f}{\partial x_i}$ if char $(K) >0$, we compute the stabilzer of $G$ and its dimension with the help of Gr\"obner bases.\\

\textbf{Problem:} Does finite determinacy of $A \in M_{r,s}$ always imply finite codimension of $\widetilde{T}_A(GA)$ if $p>0$?

We may conjecture that this not the case for arbitrary $r,s,n$.

%
%
%


Fachbereich Mathematik, Universit\"at Kaiserslautern, Erwin-Schr\"odinger Str.,
67663 Kaiserslautern, Germany

E-mail address: greuel@mathematik.uni-kl.de


\addcontentsline{toc}{section}{References}
\begin{thebibliography}{40}

\bibitem[AGV85]{AGV85} Arnol'd, V.I.; Gusein-Zade, S.M.; Varchenko, A.N.:
\emph{Singularities of differentiable maps. Volume I: The classification of critical points, caustics and wave fronts.} Monographs in Mathematics, Vol. 82. Boston-Basel-Stuttgart: Birkhäuser. X, 382 p. (1985). Zbl 0554.58001

\bibitem[BGM12]{BGM12}  Boubakri, Yousra; Greuel, Gert-Martin; Markwig, Thomas:  
\emph{Invariants of hypersurface singularities in positive characteristic}, Rev. Mat. Complut. 25, 61--85  (2012). Zbl 1279.14004

\bibitem[BGS87]{BGS87}  Buchweitz, R.-O.; Greuel, G.-M.; Schreyer, F.-O.:
\emph{Cohen-Macaulay modules on hypersurface singularities. II. }
Invent. Math. 88, 165-182 (1987). Zbl 0617.14034

\bibitem[Ca80]{Ca80} Campillo, Antonio: 
\emph{Algebroid curves in positive characteristic.} Lecture Notes in Mathematics 813. Berlin-Heidelberg-New York: Springer-Verlag, 168 p. (1980). Zbl 0451.14010
	
\bibitem[Ca83]{Ca83} Campillo, Antonio: 
\emph{Hamburger-Noether expansions over rings.} Trans. Am. Math. Soc. 279, 377-388 (1983). Zbl 0559.14020

\bibitem[CGL07]{CGL07} Campillo, Antonio; Greuel, Gert-Martin; Lossen, Christoph:
\emph{Equisingular deformations of plane curves in arbitrary characteristic.} Compos. Math. 143, 829-882 (2007). Zbl 1121.14003

\bibitem[DGPS16]{DGPS16} Decker, W.; Greuel, G.-M.; Pfister, G.; Sch{\"o}nemann, H.:
{\sc Singular} 4-1-0 -- \emph{{A} computer algebra system for polynomial computations}. http://www.singular.uni-kl.de (2016).

\bibitem[Du79]{Du79} Durfee, Alan H.:
\emph{Fifteen characterizations of rational double points and simple critical points.} Enseign. Math., II. Sér. 25, 132-163 (1979). Zbl 0418.14020

\bibitem[Gr92]{Gr92} Greuel, G.-M.;
\emph{Deformation and classification of singularities and modules. (Deformation und Klassifikation von Singularitäten und Moduln.)} Jahresbericht der DMV, Jubiläumstag., 100 Jahre DMV, Bremen/Dtschl. 1990, 177-238 (1992). Zbl 0767.32017

\bibitem[GK90]{GK90} Greuel, G.-M.; Kr\"oning, H.:
\emph{Simple singularities in positive characteristic.}
Math. Z. 203, No.2, 339-354 (1990). Zbl 0715.14001

\bibitem[GLS07]{GLS07} Greuel, Gert-Martin;  Lossen, Christoph; Shustin, Eigenii:
 \emph{Introduction to Singularities and Deformations}, Springer Monographs in Mathematics, Springer, Berlin, 2007.  Zbl 1125.32013		
	
\bibitem[GN16]{GN16} Greuel, Gert-Martin; Nguyen, Hong Duc:
 \emph{Right simple singularities in positive characteristic.} 
J. Reine Angew. Math. 712, 81-106 (2016). Zbl 1342.14006
	
\bibitem[GP16]{GP16} Greuel, Gert-Martin; Pham, Thuy Huong:
  \emph{On finite determinacy for matrices of power series.}  arXiv:1609.05133 (2016), to appear in Math. Zeitschrift.

\bibitem[GP17]{GP17} Greuel, Gert-Martin; Pham, Thuy Huong:
  \emph{Finite determinacy of matrices and ideals in arbitrary characteristic.} arXiv:1708.02442 (2017)
  
 \bibitem[GP17a]{GP17a} Greuel, Gert-Martin; Pham, Thuy Huong:
  \emph{ Algorithms for group actions in arbitrary characteristic and a problem in singularity theory}.  	arXiv:1709.08592

\bibitem[Kn87]{Kn87} Kn\"orrer, Horst:
 \emph{Cohen-Macaulay modules on hypersurface singularities. I.} 
Invent. Math. 88, 153-164 (1987). Zbl 0617.14033

\bibitem[Ng17]{Ng17} Nguyen, Hong Duc:
 \emph{Right unimodal and bimodal singularities In positive characteristic},  arXiv:1507.03554 (2017)

\end{thebibliography}
\end{document}